\documentclass[12pt]{amsart}
\usepackage{amssymb}
\usepackage{mathrsfs}
\usepackage{amsxtra}
\usepackage{graphics}
\usepackage{latexsym}
\usepackage{xcolor}
\usepackage{amsmath}
\usepackage{amssymb,amsthm,amsfonts}
\usepackage{mathtools}
\usepackage{amscd}
\usepackage{tikz-cd}
\usepackage[arrow, matrix, curve]{xy}
\usepackage{syntonly}
\usepackage{paralist}
\usepackage[plainpages=false,pdfpagelabels,
colorlinks=true,linkcolor=blue,
citecolor=blue,filecolor=blue,      
urlcolor=blue,hypertexnames=false]{hyperref}
\usepackage{cleveref}
\usepackage{quiver} 
\usepackage{stmaryrd} \SetSymbolFont{stmry}{bold}{U}{stmry}{m}{n}   
\title[]{Nonabelian Kodaira-Spencer maps}

\author[Yixuan Fu]{Yixuan Fu}
\email{fuyx23@mails.tsinghua.edu.cn}
\address{Yau Mathematical Science Center, Tsinghua University, Beijing, 100084, China}
 
\author[Mao Sheng]{Mao Sheng}
 
\email{msheng@tsinghua.edu.cn}
\address{Yau Mathematical Science Center, Tsinghua University, Beijing, 100084, China}
\address{Yanqi Lake Beijing Institute of Mathematical Sciences and Applications, Beijing, 101408, China}

\begin{document}
\theoremstyle{plain}
\newtheorem{thm}{Theorem}[section]
\newtheorem{theorem}[thm]{Theorem}
\newtheorem{lemma}[thm]{Lemma}
\newtheorem{corollary}[thm]{Corollary}
\newtheorem{proposition}[thm]{Proposition}
\newtheorem{addendum}[thm]{Addendum}
\newtheorem{variant}[thm]{Variant}
\theoremstyle{definition}
\newtheorem{lemma and definition}[thm]{Lemma and Definition}
\newtheorem{construction}[thm]{Construction}
\newtheorem{statement}[thm]{Statement}
\newtheorem{notations}[thm]{Notations}
\newtheorem{question}[thm]{Question}
\newtheorem{problem}[thm]{Problem}
\newtheorem{remark}[thm]{Remark}
\newtheorem{remarks}[thm]{Remarks}
\newtheorem{definition}[thm]{Definition}
\newtheorem{claim}[thm]{Claim}
\newtheorem{assumption}[thm]{Assumption}
\newtheorem{assumptions}[thm]{Assumptions}
\newtheorem{properties}[thm]{Properties}
\newtheorem{example}[thm]{Example}
\newtheorem{conjecture}[thm]{Conjecture}
\newtheorem{proposition and definition}[thm]{Proposition and Definition}
\numberwithin{equation}{thm}
\newcommand{\Spec}{\mathrm{Spec}}
\newcommand{\pP}{{\mathfrak p}}
\newcommand{\sA}{{\mathcal A}}
\newcommand{\sB}{{\mathcal B}}
\newcommand{\sC}{{\mathcal C}}
\newcommand{\sD}{{\mathcal D}}
\newcommand{\sE}{{\mathcal E}}
\newcommand{\sF}{{\mathcal F}}
\newcommand{\sG}{{\mathcal G}}
\newcommand{\sH}{{\mathcal H}}
\newcommand{\sI}{{\mathcal I}}
\newcommand{\sJ}{{\mathcal J}}
\newcommand{\sK}{{\mathcal K}}
\newcommand{\sL}{{\mathcal L}}
\newcommand{\sM}{{\mathcal M}}
\newcommand{\sN}{{\mathcal N}}
\newcommand{\sO}{{\mathcal O}}
\newcommand{\sP}{{\mathcal P}}
\newcommand{\sQ}{{\mathcal Q}}
\newcommand{\sR}{{\mathcal R}}
\newcommand{\sS}{{\mathcal S}}
\newcommand{\sT}{{\mathcal T}}
\newcommand{\sU}{{\mathcal U}}
\newcommand{\sV}{{\mathcal V}}
\newcommand{\sW}{{\mathcal W}}
\newcommand{\sX}{{\mathcal X}}
\newcommand{\sY}{{\mathcal Y}}
\newcommand{\sZ}{{\mathcal Z}}
\newcommand{\A}{{\mathbb A}}
\newcommand{\B}{{\mathbb B}}
\newcommand{\C}{{\mathbb C}}
\newcommand{\D}{{\mathbb D}}
\newcommand{\E}{{\mathbb E}}
\newcommand{\F}{{\mathbb F}}
\newcommand{\G}{{\mathbb G}}
\renewcommand{\H}{{\mathbb H}}
\newcommand{\I}{{\mathbb I}}
\newcommand{\J}{{\mathbb J}}
\renewcommand{\L}{{\mathbb L}}
\newcommand{\M}{{\mathbb M}}
\newcommand{\N}{{\mathbb N}}
\renewcommand{\P}{{\mathbb P}}
\newcommand{\Q}{{\mathbb Q}}
\newcommand{\Qbar}{\overline{\Q}}
\newcommand{\R}{{\mathbb R}}
\newcommand{\SSS}{{\mathbb S}}
\newcommand{\T}{{\mathbb T}}
\newcommand{\U}{{\mathbb U}}
\newcommand{\V}{{\mathbb V}}
\newcommand{\W}{{\mathbb W}}
\newcommand{\Z}{{\mathbb Z}}
\newcommand{\g}{{\gamma}}
\newcommand{\id}{{\rm id}}
\newcommand{\rk}{{\rm rank}}
\newcommand{\END}{{\mathbb E}{\rm nd}}
\newcommand{\End}{{\rm End}}
\newcommand{\Hom}{{\rm Hom}}
\newcommand{\Hg}{{\rm Hg}}
\newcommand{\tr}{{\rm tr}}
\newcommand{\Sl}{{\rm Sl}}
\newcommand{\GL}{{\rm Gl}}
\newcommand{\Cor}{{\rm Cor}}

\newcommand{\SO}{{\rm SO}}
\newcommand{\OO}{{\rm O}}
\newcommand{\SP}{{\rm SP}}
\newcommand{\Sp}{{\rm Sp}}
\newcommand{\UU}{{\rm U}}
\newcommand{\SU}{{\rm SU}}
\newcommand{\SL}{{\rm SL}}
\newcommand{\ra}{\rightarrow}
\newcommand{\la}{\leftarrow}
\newcommand{\Gal}{\mathrm{Gal}}
\newcommand{\Res}{\mathrm{Res}}
\newcommand{\Gl}{\mathrm{Gl}}
\newcommand{\Gr}{\mathrm{Gr}}
\newcommand{\Exp}{\mathrm{Exp}}
\newcommand{\Sym}{\mathrm{Sym}}
\newcommand{\Ann}{\mathrm{Ann}}
\newcommand{\GSp}{\mathrm{GSp}}
\newcommand{\Tr}{\mathrm{Tr}}
\newcommand{\HIG}{\mathrm{HIG}}
\newcommand{\MIC}{\mathrm{MIC}}
\newcommand{\NHIG}{\mathrm{NHIG}}
\newcommand{\NMIC}{\mathrm{NMIC}}
\newcommand{\FV}{\mathrm{FV}}
\newcommand{\HV}{\mathrm{HV}}
\newcommand{\Ext}{\mathrm{Ext}}
\newcommand{\bA}{\mathbf{A}}
\newcommand{\bK}{\mathbf{K}}
\newcommand{\bM}{\mathbf{M}} 
\newcommand{\bP}{\mathbf{P}}
\newcommand{\bC}{\mathbf{C}}
\newcommand{\NMF}{\mathrm{NMF}}
\newcommand{\sFV}{\mathrm{SFV}}
\newcommand{\sHV}{\mathrm{SHV}}
\newcommand{\Aut}{{\rm Aut}}
\newcommand{\Gm}{{\mathbb{G}_m}}

\begin{abstract}
We give an explicit formula of the associated graded map to the nonabelian Gauss-Manin connection with respect to the nonabelian Hodge filtration.  
\end{abstract}

\maketitle

\section{Introduction}

In Hodge theory, it is a basic and extremely useful result of P. Griffiths that the Gauss-Manin connection obeys a transversality condition with respect to the Hodge filtration, and the associated graded map is nothing but the Kodaira-Spencer map coupled with the Hodge bundles (see e.g. \cite[Proposition 10.18, Theorem 10.21]{Voisin}.) Griffiths's result is in the complex analytical setting. A purely algebraic treatment is given in the preliminary part of the celebrated work \cite{Katz} of N. Katz on $p$-curvature conjecture (see \cite[Propositions 1.4.1.6-1.4.1.7]{Katz}). What would an analogue in the nonabelian Hodge theory look like? 

In \cite{Simpson} (especially \S5, \S7-\S8), C. Simpson established such an analogue, except that there is lack of an explicit formula of the \emph{nonabelian Kodaira-Spencer map}. Let us recall briefly what Simpson has done in loc. cit. Let $k$ be an algebraically closed field. Let $G$ be a connected reductive algebraic group over $k$. Let $S$ be a smooth variety over $k$, and $\alpha: X\to S$ a smooth projective morphism over $k$. Let $f: M^{sm}_{dR}(X/S,G)\to S$ be the relative de Rham moduli space parameterizing relative Zariski dense flat connections. Note that $f$ is smooth. We summarize Simpson's results as follows.
\begin{theorem}[Simpson]\label{Simpson's result}
Assume $k=\mathbb{C}$. The morphism $f$ extends to a smooth morphism $$\tilde f: M^{sm}_{Hod}(X/S,G)\to S\times \mathbb{A}^1,$$ together with a $\mathbb{G}_m$-action on $M^{sm}_{Hod}(X/S,G)$, which covers the standard $\mathbb{G}_m$-action on $S\times \mathbb{A}^1$ and whose fiber over $S\times \{0\}$ is $g: M^{sm}_{Dol}(X/S,G)\to S$, the relative Dolbeault moduli space parameterizing relative Zariski dense Higgs bundles with rational vanishing Chern classes. The nonabelian Gauss-Manin connection on $f$, which is a sheaf morphism $f^*T_S\to T_{M^{sm}_{dR}(X/S,G)}$, extends to a $\mathbb{G}_m$-equivariant sheaf morphism 
$$
\tilde \nabla: \tilde f^*T_{S\times \mathbb{A}^1/\mathbb{A}^1}(-S\times \{0\})\to T_{M^{sm}_{Hod}(X/S,G)/\mathbb{A}^1},
$$
whose composite with the natural projection $T_{M^{sm}_{Hod}(X/S,G)/\mathbb{A}^1}\to \tilde f^*T_{S\times \mathbb{A}^1/\mathbb{A}^1}$ is the natural inclusion.
\end{theorem}
In this article, we aim to understand the residue of $\tilde \nabla$ along $S\times \{0\}$, which is the nonabelian Kodaira-Spencer map in the title. As explained in \cite[Lemma 7.1]{Simpson} and the paragraph thereafter, it is clear that it should yield a certain Higgs structure on $g$. Our main result shows that it is indeed the case, and the nonabelian Kodaira-Spencer map admits an explicit description involving the Kodaira-Spencer map of $\alpha: X\to S$.   

Set $M^{sm}_{Dol}=M^{sm}_{Dol}(X/S,G)$. Consider the universal principal Higgs $G$-bundle $(\mathcal{E},\Theta)$ over $N_{Dol}$, where $N_{Dol}$ appears in the following Cartesian diagram
 
\[
		\xymatrix{ N_{Dol}\ar[d]_{g'}\ar[r]^{\alpha'} &M^{sm}_{Dol}\ar[d]^{g} \\
			X\ar[r]^{\alpha} & S.}
		\]
We write the Higgs field in the form $\Theta: T_{N_{Dol}/M^{sm}_{Dol}}\to ad(\mathcal{E})$, where $ad(\mathcal{E})$ is the adjoint bundle of $\mathcal{E}$. Although the universal family exists only etale locally, the bundle $ad(\mathcal{E})$ is actually well-defined. Let $$\Theta^{ad}: ad(\mathcal{E})\to ad(\mathcal{E})\otimes \Omega_{N_{Dol}/M^{sm}_{Dol}}$$ be the adjoint Higgs field. It follows from the integrality of $\Theta$ that 
$$
\Theta^{ad}\circ \Theta: T_{N_{Dol}/M^{sm}_{Dol}}\to ad(\mathcal{E})\otimes \Omega_{N_{Dol}/M^{sm}_{Dol}}.
$$
is zero. Let $\Omega^*_{Hig}(ad(\mathcal{E}),\Theta^{ad})$ be the Higgs complex of $(ad(\mathcal{E}),\Theta^{ad})$ which is expressed as
$$
ad(\mathcal{E})\stackrel{\Theta^{ad}}{\longrightarrow} ad(\mathcal{E})\otimes \Omega_{N_{Dol}/M^{sm}_{Dol}} \stackrel{\Theta^{ad}}{\longrightarrow}\cdots.
$$
Then we obtain a morphism of complexes $$(T_{N_{Dol}/M^{sm}_{Dol}},0)\to \Omega^*_{Hig}(ad(\mathcal{E}),\Theta^{ad}).$$ As $T_{N_{Dol}/M^{sm}_{Dol}}\cong g'^*T_{X/S}$, Grothendieck spectral sequence gives a natural morphism
$$
R^1\alpha_*T_{X/S}\otimes g'_*\mathcal{O}_{N_{Dol}}\cong R^1\alpha_*g'_*g'^*T_{X/S}\to  g_*\mathbb{R}^1\alpha'_*\Omega^*_{Hig}(ad(\mathcal{E}),\Theta^{ad})\cong g_*T_{M^{sm}_{Dol}/S}.
$$
So we have the natural morphism $\tau:  R^1\alpha_*T_{X/S} \to g_*T_{M^{sm}_{Dol}/S}$. 

\begin{theorem}\label{main result}
Notation as in Theorem \ref{Simpson's result}. Then the nonabelian Kodaira-Spencer map is the sheaf morphism $$\theta_{KS}=\tau\circ \rho_{KS}: T_{S}\to g_*T_{M^{sm}_{Dol}(X/S,G)/S},$$
where $\tau$ is defined as above and $\rho_{KS}: T_{S}\to R^1\alpha_*T_{X/S}$ is the Kodaira-Spencer map of $\alpha$. Moreover, it satisfies the following two properties:
\begin{itemize}
    \item [(i)] $[\theta_{KS},\theta_{KS}]=0$. That is, for $a,b\in T_S$, it holds that $[\theta_{KS}(a),\theta_{KS}(b)]=0$.
    \item [(ii)] For any $t\in \mathbb{G}_m(k)$, there is a commutative diagram
\[\begin{tikzcd}
	{g^*T_S} & {t^*g^*T_S} \\
	{T_{M_{Dol}/S}} & {t^*T_{M_{Dol}/S}.}
	\arrow[from=1-1, to=1-2]
	\arrow["{t\theta_{KS}}"', from=1-1, to=2-1]
	\arrow["{t^*\theta_{KS}}", from=1-2, to=2-2]
	\arrow["{t_*}", from=2-1, to=2-2]
\end{tikzcd}\]
\end{itemize}
\end{theorem}
In \cite{ChenTing}, T. Chen studied the same problem for a local universal curve around a very general point $[s_0]\in \mathcal{M}_g, \ g\geq 2$. Namely, locally around $s_0\in S$, the Kodaria-Spencer map $\rho_{KS}$ is an isomorphism. Therefore, her result is more about a concrete description of $\tau$ in this situation. 

Motivated by the above result, we introduce the following notion.
\begin{definition}
Let $S$ be a smooth variety over $k$. A nonlinear Higgs bundle over $S$ is a pair $(f,\theta)$, where $f: X\to S$ is a smooth morphism, and $\theta: T_S\to f_*T_{X/S}$ is a sheaf morphism satisfying the integrability condition $[\theta,\theta]=0$. A morphism of nonlinear Higgs bundles over $S$ is an arrow $h: (f,\theta)\to (g,\eta)$, where $h: X\to Y$ is an $S$-morphism making the following diagram commutes:
\[
		\xymatrix{ f^*T_S\ar[d]_{\theta}\ar[r]^{=} &h^*g^*T_S\ar[d]^{h^*\eta} \\
			T_{X/S}\ar[r]^{dh} & h^*T_{Y/S}.}
		\]
A nonlinear Higgs bundle $(f,\theta)$ is said to be graded, if for any $t\in \mathbb{G}_m(k)$, there is an isomorphism $(f,\theta)\cong (f,t\theta)$ of nonlinear Higgs bundles.  
\end{definition}
Thus our result says that $(M^{sm}_{Dol}(X/S,G)\to S, \theta_{KS})$ is a graded nonlinear Higgs bundle. For a Higgs bundle $(\mathcal{E},\theta)$ over $S$, one may attach the geometric vector bundle $\pi: E=\mathbf{Spec}(S(\mathcal{E^*}))\to S$, and by the natural isomorphism $T_{E/S}\cong \pi^*\mathcal{E}$, the Higgs field $\theta: T_S\to \mathcal{E}nd(\mathcal{E})\cong \mathcal{E}\otimes \mathcal{E}^{*}$ gives rise to a Higgs field 
$$
T_S\to \pi_*T_{E/S}\cong \mathcal{E}\otimes \pi_*\mathcal{O}_E\cong \mathcal{E}\otimes S(\mathcal{E^*})
$$
in the new sense. Thus the category of Higgs bundles is naturally enlarged into the category of nonlinear Higgs bundles, which contains nonabelian Kodaira-Spencer maps.

{\bf Acknowledgements.} The authors would like to thank Daniel Litt and Carlos Simpson for their interest in the work. This work is partially supported by the Chinese Academy of Sciences Project for Young Scientists in Basic Research (Grant No. YSBR-032).

\section{Differential calculus for $\lambda$-connections}
 
In this section, we work over a base ring $k$ with additional condition that $2$ is  regular in $k$. Our goal is to extend the differential calculus (see e.g Berthelot-Ogus \S2)
to $\lambda$-connections. All rings are assumed to be commutative and unital.
    \begin{lemma}\label{ringmorinterpolation}
        Let $f_1,\dots,f_n$ be ring homomorphisms from $A$ to $B$, such that the image of $f_i - f_j$ are contained in a square zero ideal in $B$. Then for any collection of elements $\lambda_1, \dots, \lambda_n \in B$ such that $\sum_i \lambda_i = 1$, the map $\sum_i \lambda_i f_i$ is still a ring homomorphism.
    \end{lemma}
    \begin{proof}
        We only need to verify the multiplication. Indeed,
        \begin{align*}
            (\sum_{i}\lambda_i f_i(x)) (\sum_i \lambda_i f_i(y)) - (\sum_i \lambda_i) (\sum_i \lambda_i f_i(xy))
            =\sum_{i,j}\lambda_i \lambda_j(f_i(x)-f_j(x))(f_i(y)-f_j(y)).
        \end{align*}
    \end{proof}
    Let $k \to A$ be a smooth morphism, and $\lambda \in k$. We also use $X,S$ to indicate $\operatorname{Spec}A,\operatorname{Spec}k$ resp..
    
    \begin{definition}
     A $\lambda$-connection on an $A$-module $E$ is a $k$-linear homomorphism
    \[
    \nabla \colon E \to E \otimes_A \Omega_{A/k}^1
    \]
    such that
    \[
    \nabla (fs) = s \otimes \lambda df + f\nabla s, \ f \in A, s \in E.
    \]  
    It is said to be integrable if the curvature 
$$
\Theta_{\nabla}=\nabla^1\circ \nabla: E\to E\otimes_A\Omega^2_{A/k} 
$$
vanishes, where $\nabla^1: E\otimes \Omega^1_{A/k}\to E\otimes \Omega^2_{A/k}$ is defined by 
$$
s\otimes \omega\mapsto -s\otimes \lambda d\omega+\nabla(s)\wedge \omega, \ s\in E, \omega\in \Omega^1_{A/k}.
$$

If $x_1,\dots,x_n$ are local coordinates of the smooth $k$-algebra $A$, we use $\nabla_{\frac{\partial}{\partial x_i}}$ or $\nabla_i$ to indicate the morphism
\[\begin{tikzcd}
	E & {E \otimes \Omega^1_{A/k}} & E.
	\arrow["\nabla", from=1-1, to=1-2]
	\arrow["{\langle \cdot, \frac{\partial}{\partial x_i} \rangle}", from=1-2, to=1-3]
\end{tikzcd}\]
    \end{definition}

    Let $B$ be another $k$-algebra, and $f_0,f_1$ be two homomorphisms from $A$ to $B$ such that $f_1 - f_0$ is contained in a square-zero ideal of $B$. A typical example is $B = P^1$, the coordinate ring of the first order neighborhood $D(1)$ of $X$ in its diagonal $X\times_SX$, and the two morphisms are two projections $d_0,d_1: D(1)\to X$. From Lemma \ref{ringmorinterpolation} we know $f_\lambda := (1 - \lambda)f_0 + \lambda f_1$ is a ring homomorphism. Notice that 
    \[
    f_0 = (1- \lambda)f_\lambda + \lambda f_{\lambda - 1}.
    \]
    \begin{proposition}
        A $\lambda$-connection on $E$ is equivalent to an isomorphism $\epsilon \colon d_\lambda^*E \simeq d_0^*E$ which reduces to the identity modulo diagonal.
    \end{proposition}
    \begin{proof}
        Given $\nabla$ we define a homomorphism $\theta \colon E \to d_{\lambda*}d_0^*E$ by
        \[
        \theta (x) = \nabla x +  x \otimes \mathbf{1}.
        \]
        We need to verify $\theta$ is linear, namely
        \[
        \theta (fx) = d_\lambda(f)\theta(x).
        \]
        Indeed, 
        \[
        \theta(fx) = \nabla (fx) + fx \otimes \mathbf{1} = \lambda x \otimes df + f\nabla x + fx \otimes 1,
        \]
        \[
        d_\lambda (f)\theta(x) = (f \otimes 1 + \lambda df)(\nabla x + x \otimes \mathbf{1}).
        \]
        Notice that $df \nabla x = 0$, the result follows. 

        By extension of scalars, we get $\epsilon \colon d_\lambda^*E\to d_0^*E$. It is an isomorphism since its reduction to diagonal is the identity map.

        Conversely, from $\epsilon$ we get the corrsesponding restriction of scalars $\theta$, and define $\nabla(x) = \theta(x) - x \otimes \mathbf{1} $.
    \end{proof}
\begin{example}
   A typical example of $\lambda$-connection is $S = \mathbb{A}^1$ and $X = X_0 \times \mathbb{A}^1$, and $\lambda \in \Gamma(\mathcal{O}_{\mathbb{A}^1})$ is the standard coordinate. In this case, a $\lambda$-connection on a vector bundle $E$ over $X_0 \times \mathbb{A}^1$ is given by an isomorphism $p_t^*E \simeq p_0^*E$, where $p_t$ is given by the ring homomorphism
    \begin{align*}
        A[t] \to P^1[t], \qquad a \mapsto (1-t)d^*_0a + td_1^*a, \quad t \mapsto t.
    \end{align*}
\end{example}

    \begin{corollary}
Let $\nabla$ be a $\lambda$-connection on an $A$-module $E$. If $f_0,f_1$ are two $k$-algebra morphisms from $A$ to $B$ such that $f_1 - f_0$ is contained in a square-zero ideal in $B$, then there is an isomorphism of $B$-modules $\epsilon \colon f_\lambda^*E \simeq f_0^*E$ with in local coordinates $x_i$ of the smooth $k$-algebra $A$
        \[
        s \otimes 1 \mapsto s \otimes 1 +\sum \nabla_i s \otimes (f_1^*x_i - f_0^*x_i).
        \]
        
    \end{corollary}
    \begin{proposition}\label{pullbackconnection}
        If $(E,\nabla)$ is integrable, then $\epsilon: f_\lambda^*(E,\nabla) \simeq f_0^*(E,\nabla)$. That is, the isomorphism $\epsilon$ constructed above is not only an isomorphism of sheaves, but also preserves the pullback connections.
    \end{proposition}
    \begin{proof}
           The pullback connection is given by
        \[
        f^* \nabla (f^*s) = \sum_i f^* \nabla_i s \otimes d f^*x_i.
        \]
        From the expression of $\epsilon$ we also get
        \[
        \epsilon(\nabla_i s \otimes 1) = \nabla_is \otimes 1 + \nabla_j\nabla_is \otimes (f_1^*x_j - f_0^*x_j).
        \]
        Therefore we obtain
        \begin{align}\label{pullbackconnection1}
            \epsilon (f_\lambda^*\nabla(s \otimes 1)) = \nabla_i s \otimes 1 \otimes df_\lambda^*x_i +\nabla_j \nabla_i s \otimes (f_1^*x_j - f_0^*x_j) \otimes df_\lambda^*x_i.
        \end{align}
       On the other hand, 
       \begin{align*}
           f_0^*\nabla(s \otimes 1) = \nabla_i s \otimes 1 \otimes df_0^*x_i,
       \end{align*}
       \begin{align}\label{pullbackconnection2}
           f_0^*\nabla (\nabla_is \otimes (f_1^*x_i - f_0^*x_i))
           = \lambda \nabla_is \otimes d(f_1^*x_i - f_0^*x_i) + f_0^*\nabla(\nabla_is) \otimes (f_1^*x_i - f_0^*x_i),
       \end{align}
       where 
       \begin{align*}
           f_0^*\nabla(\nabla_is) = \nabla_j\nabla_is \otimes 1 \otimes df_0^*x_j.
       \end{align*}
        Collecting the first order term of $\nabla$ in (\ref{pullbackconnection2}), we get
        \[
        \nabla_is \otimes 1 \otimes f_0^*x_i + \lambda \nabla_i s \otimes 1 \otimes d(f_1^*x_i-f_0^*x_i) = \nabla_is \otimes 1 \otimes f_\lambda^*x_i.
        \]
        Comparing the second order term of $\nabla$ in (\ref{pullbackconnection1}) and (\ref{pullbackconnection2}) we then need to show
        \[
        \nabla_j\nabla_is \otimes (f_1^*x_j - f_0^*x_j) \otimes df_\lambda^*x_i = \nabla_j\nabla_is \otimes (f_1^*x_i -f_0^*x_i) \otimes df_0^*x_j.
        \]
        Since $\nabla_i \nabla_j = \nabla_j\nabla_i$, the difference is
        \[
        \nabla_j\nabla_is \otimes (f_1^*x_j-f_0^*x_j)(df\lambda^*x_i - df_0^*x_i)
        = \lambda  \nabla_j\nabla_is \otimes (f_1^*x_j-f_0^*x_j)(d(f_1^*x_i - f_0^*x_i))
        \]
        Notice that since $(f_1^*x_j-f_0^*x_j)(f_1^*x_i - f_0^*x_i) = 0$, we have
        \[
        (f_1^*x_j-f_0^*x_j)(d(f_1^*x_i - f_0^*x_i)) + (f_1^*x_i-f_0^*x_i)(d(f_1^*x_j - f_0^*x_j)) = 0.
        \]
        The result follows.
    \end{proof}
       \begin{proposition}\label{planergeometry}
        Let $f_0,f_1,f_2,f_{01},f_{12},f_{20}$ be ring morphisms from $A$ to $B$ such that their pairwise differences are contained in some square zero ideal of $B$, and
        \begin{align*}
            f_1 = (1-\lambda) f_0 + \lambda f_{01},\\
            f_2 = (1- \lambda)f_1 + \lambda f_{12},\\
            f_0 = (1-\lambda)f_2 + \lambda f_{20}.
        \end{align*}
        If moreover
        \[
        f_0 + f_1 + f_2 = f_{01} + f_{12}+ f_{20},
        \]
        then the composition
        \[\begin{tikzcd}
	{f_0^*E} & {f_1^*E} & {f_{2}^*E} & {f_0^*E}
	\arrow["{\epsilon_{f_0,f_{01}}}", from=1-1, to=1-2]
	\arrow["{\epsilon_{f_1,f_{12}}}", from=1-2, to=1-3]
	\arrow["{\epsilon_{f_2,f_{20}}}", from=1-3, to=1-4]
\end{tikzcd}\]
is equal to the identity.
    \end{proposition}
    \begin{proof}
        From the explicit expression of $\epsilon$ we see the composition is 
        \[
        s \otimes 1 \mapsto s \otimes 1 + \sum\nabla_i s \otimes (f_0+f_1+f_2-f_{01}-f_{12}-f_{02})^*x_i.
        \]
    \end{proof}

Let $B$ be a commutative $A$-algebra. A $\lambda$-transversal distribution on $B$ over $A$ is an isomorphism of $P_A^1$-algebras
\[
\epsilon \colon d_\lambda^*B \simeq d_0^*B.
\]
which reduces to identity modulo $\Omega_{A/k}^1$.
\begin{proposition}
    A $\lambda$-transversal distribution on $B$ is equivalent to a left $\lambda$-splitting of the exact sequence 
    \[\begin{tikzcd}
	{B \otimes_A\Omega_{A/k}} & {\Omega_{B/k}} & {\Omega_{B/A}} & 0.
	\arrow[from=1-1, to=1-2]
	\arrow[from=1-2, to=1-3]
	\arrow[from=1-3, to=1-4]
\end{tikzcd}\]
\end{proposition}
\begin{proof}
Regard $B$ as an $A$-module we see $B$ admits a $\lambda$-connection $\nabla$. The restriction of scalars $\theta \colon B \to d_{\lambda*} d_0^*B$ preserves multiplication, namely
\[
\nabla (xy) + xy \otimes 1 = (\nabla x + x \otimes 1)(\nabla y + y \otimes 1).
\]
Therefore, the map $\nabla \colon B \to B \otimes \Omega_{A/k}^1$ satisfies
\[
\nabla(xy) = x \nabla y + y \nabla x,
\]
namely $\nabla \in \operatorname{Der}_{k}(B,B \otimes \Omega_{A}^1)$. This gives rise to a morphism $\Omega_B \to B \otimes \Omega_A$. For $x,y \in A$ we have
\[
\nabla (xy) = \lambda  \otimes d(xy).
\]
Thus $\nabla$ is a $\lambda$-splitting. The converse procedure is also clear.
\end{proof}

Everything we discussed above can be globalized to schemes. Let $S$ be a separated smooth $k$-scheme. We denote $D_S(1)$ for the first order infinitesimal neighborhood of  the diagonal $S \to S \times_k S$. We have a truncated simplicial diagram of schemes
\[\begin{tikzcd}
	{D(1)} & S.
	\arrow["{p_0}", shift left=3, from=1-1, to=1-2]
	\arrow["{p_1}"', shift right=3, from=1-1, to=1-2]
	\arrow["\Delta"{description}, from=1-2, to=1-1]
\end{tikzcd}\]

 If $E$ is a quasi-coherent $\mathcal{O}_S$-module, a $\lambda$-connection on $E$ relative to $k$ is an isomorphism of $D_S(1)$-modules $\epsilon \colon p_\lambda^* E \simeq p_0^*E$ that is $\mathrm{id}_E$ after restriction to the diagonal. If $X$ is a scheme over $S$, a $\lambda$-transversal distribution of $X/S$ is an isomorphism of
 $D_S(1)$-schemes $\epsilon \colon p_0^* X \simeq p_\lambda^*X$ that is $\mathrm{id}_X$ after restriction to the diagonal. Equivalently, it is a left $\lambda$-splitting of the short exact sequence
\[\begin{tikzcd}
	{f^* \Omega_{S/k}} & {\Omega_{X/k}} & {\Omega_{X/S}} & 0.
	\arrow[from=1-1, to=1-2]
	\arrow[from=1-2, to=1-3]
	\arrow[from=1-3, to=1-4]
\end{tikzcd}\]
\begin{remark}
    If $S$ is smooth over $k$, then $\Omega_{S/k}$ is locally free, and thus a $\lambda$-splitting of the above sequence is equivalent to a right $\lambda$-splitting of the exact sequence
\[\begin{tikzcd}
	{T_{X/S}} & {T_{X/k}} & {f^*T_{S/k}} & 0.
	\arrow[from=1-1, to=1-2]
	\arrow[from=1-2, to=1-3]
	\arrow[from=1-3, to=1-4]
\end{tikzcd}\]
\end{remark}

We define $\mathbf{T}_{S/k} := \operatorname{Spec}(\operatorname{Sym} \Omega_{S/k})$. From the universal property of symmetric algebra and differential module, we know
\[
\mathbf{T}_{S/k}(Z) = S(Z[\epsilon]).
\]
A $
\lambda$-transversal distribution on $X$ gives a section of the morphism of tangent schemes
$f^*\mathbf{T}_{S/k} \to \mathbf{T}_{X/k}$. In terms of functor of points, it is described as follows: given a $Z$-point of $f^*\mathbf{T}_S$ 
which may be described via the following commutative diagram with solid arrows:

\[\begin{tikzcd}
	Z & X \\
	{Z[\epsilon]} & S
	\arrow["t", from=1-1, to=1-2]
	\arrow[hook, from=1-1, to=2-1]
	\arrow[from=1-2, to=2-2]
	\arrow[dashed, from=2-1, to=1-2]
	\arrow["\delta", from=2-1, to=2-2]
\end{tikzcd}\]
we want to find a point in $\mathbf{T}_{X/k}(Z)$, that is, a dashed arrow filling the diagram. 
Let 
\[
\iota \colon Z[\epsilon] \to Z \to Z[\epsilon] \stackrel{\delta}\to S
\]
be the ``constant'' map induced by $\delta$. Denote $\delta_\lambda = (1- \lambda) \iota + \lambda \delta$, then we have an isomorphism $\epsilon \colon    \iota^*X  \simeq \delta_\lambda^*X$ . 
Now $t$ gives a section $\alpha \colon Z[\epsilon] \to \iota^*X$, then composing the isomorphism $\iota^*X \simeq \delta_\lambda^*X$ we get a map $Z[\epsilon] \to \delta_\lambda^*X$, which is a $\lambda$-lifting of $\delta$.
We will call this a first order horizontal lift.

\begin{proposition}
    The first order horizontal lift is given by a $\lambda$-splitting of the exact sequence of differential forms. 
\end{proposition}
\begin{proof}
We may reduce everything to the affine case, and hence assume that $k = \operatorname{Spec} k$, $S = \operatorname{Spec} A$, $X = \operatorname{Spec} B$ and $Z = \operatorname{Spec} C$. Firstly, we shall note that the isomorphism
\[
\mathrm{Hom}_k(\operatorname{Sym}_A \Omega_{A/k},C) \simeq \mathrm{Hom}_k(A,C[\epsilon])
\]
should be interpreted as follows. Both sides are naturally isomorphic to the data of pairs $(g,D)$ where $g$ is  a $k$-algebra morphism $g \colon A \to C$ and $D$ is a derivation of $A$ valued in $C$, whose $A$-module structure is given by $g$. Indeed, homomorphisms from $A$ to $C[\epsilon]$ are of forms $g + \epsilon D$, and a homomorphism from $\operatorname{Sym}_A(\Omega_{A/k})$ to $C$ are determined by $(g,D)$ where $g$ is the map on degree $0$ parts, and $D$ is the map on degree $1$ parts, mapping $dx$ to $Dx$.

With this in mind, given a precrystal structure $\nabla \colon B \to B \otimes_A \Omega_{A/k}$, the morphism
$\Omega_{B/k}\to B \otimes \Omega_{A/k}$ is given by $db \mapsto \nabla b$.
Given a commutative diagram with solid arrows
\[\begin{tikzcd}
	B & C \\
	A & {C[\epsilon]}
	\arrow["t", from=1-1, to=1-2]
	\arrow[dashed, from=1-1, to=2-2]
	\arrow["f", from=2-1, to=1-1]
	\arrow["{(g,D)}", from=2-1, to=2-2]
	\arrow[from=2-2, to=1-2]
\end{tikzcd}\]
which is a $C$-point in $B \otimes_A\operatorname{Sym}_A(\Omega_{A})$, we seek for a dashed diagram filler which means a $C$-point of $\operatorname{Sym}_B(\Omega_{B})$. Suppose $\nabla b = \sum b_i \otimes da_i$, then under the morphism $\Omega_B \to B \otimes \Omega_A \to  C$, we have 
\begin{equation}\label{splitting}
db \longmapsto \sum b_i \otimes da_i \longmapsto \sum b_i Da_i.
\end{equation}

On the other hand, denote $f_0 = (g,0)$ and $f_1 = (g,D)$, then we need to understand the composed morphism
\[\begin{tikzcd}
	B & {f_\lambda^*B} & {f_0^*B} \\
	&& {C[\epsilon]}
	\arrow[from=1-1, to=1-2]
	\arrow["\simeq", from=1-2, to=1-3]
	\arrow[from=1-2, to=2-3]
	\arrow["{t[\epsilon]}", from=1-3, to=2-3]
\end{tikzcd}\]
Under $d_\lambda^* B \simeq d_0^*B$, we have
\[
    \mathbf{1} \otimes b \longmapsto b \otimes \mathbf{1} + \sum b_i \otimes (1 \otimes a_i - a_i \otimes 1).
\]
After pulling this back along $(f_1,f_0)$, we obtain
\begin{equation}\label{firstorderhorizontallift}
    1 \otimes b\longmapsto b \otimes 1 + \sum b_i \otimes (f_1(a_i)-f_0(a_i)) = b + \sum b_iDa_i \epsilon.
\end{equation}
Comparing (\ref{splitting}) with (\ref{firstorderhorizontallift}), we see that they indeed coincide.
\end{proof}

Let $S$ be a smooth scheme over $k$, and $f \colon X \to S\times \mathbb{A}^1$ be a finitely presented morphism, and we denote $X_\lambda$ to be the inverse image of $S \times \{ \lambda \}$ under $f$. Let $X$ be equipped with a $\mathbb{G}_m$-action such that $f \colon X \to S \times \mathbb{A}^1$ is $\mathbb{G}_m$-equivariant, where $\mathbb{G}_m$ acts on $S \times \mathbb{A}^1$ via the action of $\mathbb{G}_m$ on $\mathbb{A}^1$. Let $\nabla_\lambda$ be a $\lambda$-connection on $X$ over $S \times \mathbb{A}^1$, namely a $\lambda$-splitting of
\[
T_{X/\mathbb{A}^1} \to f^*T_{S \times \mathbb{A}^1/\mathbb{A}^1}.
\]
Notice that these sheaves are $\mathbb{G}_m$-equivariant sheaves on $X$ by the following proposition.
\begin{proposition}
    Let $f \colon X \to Y$ be a $G$-equivairiant morphism, then $\Omega_{X/Y}$ is a $G$-equivariant sheaf on $X$. If $G$ is flat, and $f$ is of finite presentation, then $T_{X/Y}$ is also a $G$-equivariant sheaf.
\end{proposition}
\begin{proof}
Denote $\mu_X \colon G \times X \to X$ to be the action.
From the Cartesian diagram
\[\begin{tikzcd}
	{G \times X} & {G \times X} & X \\
	{G \times Y} & {G \times Y} & Y
	\arrow["{(pr_G,\mu_X)}", from=1-1, to=1-2]
	\arrow[from=1-1, to=2-1]
	\arrow["\ulcorner"{anchor=center, pos=0.125}, draw=none, from=1-1, to=2-2]
	\arrow["{pr_X}", from=1-2, to=1-3]
	\arrow[from=1-2, to=2-2]
	\arrow["\ulcorner"{anchor=center, pos=0.125}, draw=none, from=1-2, to=2-3]
	\arrow["f", from=1-3, to=2-3]
	\arrow["{(pr_G,\mu_Y)}", from=2-1, to=2-2]
	\arrow["{pr_Y}", from=2-2, to=2-3]
\end{tikzcd}\]
\end{proof}
and taking the outer rectangle, we get the desired isomorphism
\[
pr_X^* \Omega_{X/Y} \simeq \mu_X^*\Omega_{X/Y}.
\]
If $f$ is finitely presented, then $\Omega_{X/Y}$ is of finite presentation, and thus taking dual commutes with flat base change.
\begin{definition}
    We say $\nabla_\lambda$ is $\mathbb{G}_m$-equivariant if the following diagram commutes
\[\begin{tikzcd}
	{pr_X^*f^*T_{S \times \mathbb{A}^1}/\mathbb{A}^1} & {\mu_X^*f^*T_{S \times \mathbb{A}^1}/\mathbb{A}^1} \\
	{pr_X^*T_{X/\mathbb{A}^1}} & {\mu_X^*T_{X/\mathbb{A}^1}}
	\arrow["\simeq", from=1-1, to=1-2]
	\arrow["{t\cdot pr_X^* \nabla_\lambda}"', from=1-1, to=2-1]
	\arrow["{\mu_X^* \nabla_\lambda}", from=1-2, to=2-2]
	\arrow["\simeq", from=2-1, to=2-2]
\end{tikzcd}\]
where $t$ is the coordinate function on $\mathbb{G}_m$.
\end{definition}

\begin{definition}
    The residue of a $\lambda$-connection $\nabla_\lambda$ on $X$ is the base change of $\nabla_t$ along $S \times \{0\} \to S \times \mathbb{A}^1$. We denote it by $\nabla_0$, and this is a Higgs field on $X_0/S$.
\end{definition}
\begin{proposition}\label{gradedness}
    The residue of a $\mathbb{G}_m$-equivariant $\lambda$-connection is a graded Higgs field.
\end{proposition}
\begin{proof}
    Take the fiber at zero we get the diagram
\[\begin{tikzcd}
	{pr_{X_0}^*f_0^*T_{S}} & {\mu_{X_0}^*f_0^*T_{S}} \\
	{pr_{X_0}^*T_{X_0}} & {\mu_{X_0}^*T_{X_0}}.
	\arrow["\simeq", from=1-1, to=1-2]
	\arrow["{t\nabla_0}"', from=1-1, to=2-1]
	\arrow["{\mu_X^*\nabla_0}", from=1-2, to=2-2]
	\arrow["\simeq", from=2-1, to=2-2]
\end{tikzcd}\]
and the conclusion follows. This is somewhat stronger than gradedness, as the isomorphism is actually given by the action of $\mathbb{G}_m$.
\end{proof}
To justify the name residue, we need the following.
\begin{proposition}
    Suppose $f \colon X \to S \times \mathbb{A}^1$ is flat. Then there is a bijective correspondence between $t$-connections on $X/S \times \mathbb{A}^1$ and maps $\tilde \nabla \colon f^*T_{S \times \mathbb{A}^1/\mathbb{A}^1}(-D) \to T_{X/\mathbb{A}^1}$ such that the composition 
\[\begin{tikzcd}
	{f^*T_{S \times \mathbb{A}^1/\mathbb{A}^1}(-D)} & {T_{X/\mathbb{A}^1}} & {f^*T_{S \times \mathbb{A}^1/\mathbb{A}^1}}
	\arrow["{\tilde \nabla}", from=1-1, to=1-2]
	\arrow[from=1-2, to=1-3]
\end{tikzcd}\]
is the canonical inclusion. Here $D$ is the effective Cartier divisor $S \times \{0\} \subset S \times \mathbb{A}^1$.
\begin{proof}
    Given $\nabla_t$, we simply let $\tilde \nabla$ to be the restriction of $\nabla_t$ to the subsheaf $f^*T_{S \times \mathbb{A}^1/\mathbb{A}^1}(-D)$. We do not need flatnees for this side.
    
     Conversely, given $\tilde \nabla$, we define$
    \nabla_t(X) := t^{-1}\tilde \nabla(tX).
    $
    Notice that it is possible to divide by $t$ since $\tilde \nabla(tX)$ is zero on $X_0$, and $X \to \mathbb{A}^1$ is flat.
\end{proof}
\end{proposition}

\begin{remark}
    From the construction we see $\nabla_t$ is $\mathbb{G}_m$-equivariant if and only if the corresponding $\tilde \nabla  \colon f^*T_{S \times \mathbb{A}^1/\mathbb{A}^1}(-D) \to T_{X/\mathbb{A}^1}$ is a morphism of $\mathbb{G}_m$-equivariant sheaves.
\end{remark}
\begin{remark}
    If $\nabla_t$ corresponds to $\tilde \nabla$, we define $\mathrm{res}_{t = 0} \tilde \nabla := \nabla_0$.
\end{remark}

\begin{definition}
    Let $f \colon X \to S \times\mathbb{A}^1$ be as above, and $\nabla$ be a connection on $X_1$. It is said to be \textit{Griffiths transversal} to $f$ if there exists a $\mathbb{G}_m$-equivariant $t$-connection $\nabla_t$  on $X$ such that $\nabla_1 = \nabla$. In this case, we say $\nabla$ is a \textit{filtered connection}, and the residue $\nabla_0$ is said to be the associated graded Higgs field to $\nabla$.
\end{definition}

\section{Non-abelian $\lambda$-connections}
Let $S$ be a smooth $k$-scheme, and let $X \to S$ be a smooth projective morphism. Consider the Hodge moduli stack
    \[
    \mathcal{M}_{Hod} (T) = \{\lambda \in \Gamma(\mathcal{O}_T), (E,\nabla) \in \mathrm{Conn}_\lambda (X_T/T) \}.
    \]
    To get good moduli space, we need to consider the open substack consisting of semistable integrable $\lambda$-connections with vanishing rational Chern classes. For the construction of nonabelian $\lambda$-connection, this does not matter. 
    
    There is a natural morphism $\mathcal{M}_{Hod} \to \mathbb{A}_S^1$ given by projecting to $\lambda$. We are going to construct the nonabelian $\lambda$-connection on $\mathcal{M}_{Hod}$ over $\mathbb{A}^1$. Recall the diagram
\[\begin{tikzcd}
	& {D(1) \times \mathbb{A}^1} \\
	{\mathcal{M}_{Hod}} & {S \times \mathbb{A}^1}.
	\arrow["{p_{\lambda}}"', shift right=2, from=1-2, to=2-2]
	\arrow["{p_0}", shift left=2, from=1-2, to=2-2]
	\arrow[from=2-1, to=2-2]
\end{tikzcd}\]    
    We need to specify an isomorphism of stacks
    \[
    p_0^* \mathcal{M}_{Hod} \simeq p_{\lambda}^*\mathcal{M}_{Hod},
    \]
    which shall be called the non-abelian $\lambda$-connection.
    This means that, given scheme $T$ and $\lambda \in \Gamma(\mathcal{O}_T)$, with morphisms $f_0,f_1$ to $S$ which differ by a square zero ideal, we need to construct an equivalence
    \[
    \mathrm{Conn}_\lambda(f_\lambda^*X/T) \simeq \mathrm{Conn}_\lambda(f_0^*X/T).
    \]
    Let $T_0$ be the pullback
\[\begin{tikzcd}
	{T_0} & S \\
	T & {D(1)}
	\arrow[from=1-1, to=1-2]
	\arrow[from=1-1, to=2-1]
	\arrow[hook, from=1-2, to=2-2]
	\arrow["{(f_0,f_1)}", from=2-1, to=2-2]
\end{tikzcd}\]
and $X_0$ be the pullback of $X$ to $T_0$. Then $f_i^*X$ are both smooth liftings of $X_0$ to $T$, hence are locally isomorphic.
\[\begin{tikzcd}
	{X_0} & {f_i^*X} && X \\
	{T_0} & T & {D(1)} & S
	\arrow[hook, from=1-1, to=1-2]
	\arrow[from=1-1, to=2-1]
	\arrow[from=1-2, to=2-2]
	\arrow[from=1-4, to=2-4]
	\arrow[hook, from=2-1, to=2-2]
	\arrow["{(f_0,f_1)}", from=2-2, to=2-3]
	\arrow["{d_1}"', shift right, from=2-3, to=2-4]
	\arrow["{d_0}", shift left, from=2-3, to=2-4]
\end{tikzcd}\]
Take affine open covers $U_\alpha$ of $X_0$ so that we have isomorphisms
\[
f_\alpha \colon (f_1^*X)_\alpha \simeq (f_0^*X)_\alpha.
\]
Take $f_{\alpha,\lambda}(x) = (1- \lambda) x \otimes 1 + \lambda f_\alpha(x)$, then we have
\[
f_{\alpha,\lambda} \colon (f_\lambda^*X)_\alpha \simeq (f_0^*X)_\alpha.
\]

Given a $\lambda$-connection $(E,\nabla)$ on $f_0^*X$, we denote
$E_\alpha = f_{\alpha,\lambda}^*(E,\nabla)$. These are locally connections on $f_\lambda^*X$. To glue them together we need isomorphisms $E_\alpha \simeq E_\beta$ satisfying cocycle condition. Indeed, take
\[
f_{\alpha,\beta}(x) = (1-\lambda) x\otimes 1 + f_\beta(x) - (1-\lambda) f_\alpha(x),
\]
then
\[
f_{\beta,\lambda} = (1-\lambda)f_{\alpha,\lambda} + \lambda f_{\alpha\beta}.
\]
From Proposition \ref{planergeometry} we get the desired glueing data.

We give another interpretation of this equivalence. Since $f_{\alpha,\lambda}-f_{\beta,\lambda} = \lambda(f_\alpha -f_\beta)$, we obtain 
\[
f_\lambda^*X = (1-\lambda) f_0^*X + \lambda f_1^*X.
\]
\begin{proposition}\label{invariancefirstorder}
    Let $X/S_0$ be a smooth scheme, $S_0 \hookrightarrow S$ be a square-zero thickening. Suppose $X_0,X_1$ be two smooth liftings of $X$ to $S$, and $\lambda \in \Gamma(\mathcal{O}_S)$. Since the liftings of $X$ to $S$ is a torsor under $H^1(T_{X/S_0}\otimes \mathcal{I}_{S_0})$, it makes sense to define
    \[
    X_\lambda = (1-\lambda)X_0 + \lambda X_1.
    \]
    Then there is an equivalence $\mathrm{Conn}_\lambda(X_0/S) \simeq \mathrm{Conn}_\lambda(X_\lambda/S)$.
\end{proposition}
\begin{proof}
    Suppose over an open cover $U_\alpha$, we have an isomorphism
    $\alpha \colon X_1|U_\alpha \simeq X_0|U_\alpha$. Then on $U_{\alpha\beta}$ we have a derivation $D_{\alpha\beta} \in \mathrm{Der}_{S_0}(\mathcal{O}_X, \mathcal{I})$ given by $\alpha \beta^{-1} - 1$. Then $X_1 - X_0$ is represented by the cocycle $D_{\alpha\beta}$, and $X_\lambda$ is represented by $\lambda D_{\alpha\beta}$. That is, $X_\lambda$ is given by the twisting $1+\lambda D_{\alpha\beta}$.

    For $(E,\nabla) \in \mathrm{Conn}_\lambda(X_0)$, since $\alpha\beta^{-1}$ and $1$ differs by a square zero ideal, we know by Proposition \ref{pullbackconnection} that there is an isomorphsim 
    \[
    \epsilon_{\alpha\beta} \colon (E,\nabla) \simeq (1 - \lambda +\lambda \alpha\beta^{-1})^* (E,\nabla) = (1 + \lambda D_{\alpha\beta})^*(E,\nabla)
    \]
    on $U_{\alpha\beta}$. Again from Proposition \ref{planergeometry} we have
    \[
    (1+\lambda D_{\beta\gamma})^*(\epsilon_{\beta\gamma}) \circ \epsilon_{\alpha\beta} = \epsilon_{\alpha\gamma}.
    \]
    Therefore, via $\epsilon_{\alpha\beta}$ we can glue $(E,\nabla)_\alpha$ to a $\lambda$-connection on $X_\lambda$.
\end{proof}
\section{Hodge stacks}
From now on, we shall work over $k=\mathbb{C}$. Hodge stack, as introduced by Simpson, is a $\mathbb{G}_m$-equivariant deformation of the de Rham stack to the Dolbeault stack. In this section, we shall review its definition and basic properties. Our goal is to get nonabelian 1/0-connection introduced in last section by exploiting the presentation of the deRham/Dolbeault stack.
\begin{remark}
We expect a similar explaination of the nonabelian $\lambda$ connection via the Hodge stack. However, this seems to be not very easy in the classical context. Rather, one should appeal to the more flexible language of animated rings.
\end{remark}
  
Let $X$ be a finite type smooth scheme over $\operatorname{Spec}\mathbb{C}$. Our stacks will be stacks on the big etale toplogy of $X$.

The de Rham stack $X_{dR}$ is defined by the groupoid object in $
\operatorname{Shv}((Sch/\mathbb{C})_{et})$
\[\begin{tikzcd}
	{(X \times X)^{\wedge}} & X
	\arrow[shift right, from=1-1, to=1-2]
	\arrow[shift left, from=1-1, to=1-2]
\end{tikzcd}\]
where $(X \times X)^\wedge$ is the formal completion of $X \times X$ along the diagonal, with the two morphisms induced by natural projections, and the composition is induced by $p_{02} \colon X \times X \times X \to X$. From the infinitesimal lifting criterion of smoothness, it is not hard to see $X_{dR}(S) = X(S_{red})$. 
\begin{example}
    Vector bundles on $X_{dR}$ are the same as vector bundles on $X$ with integrable connections. This follows from the groupoid presentation of $X_{dR}$, and identification of integrable connections with stratified vector bundles. Therefore, an integrable connection on a vector bundle $E$ over $E$ is the same as an isomorphism $E \simeq \tilde E \times_{X_{dR}} X$ for a vector bundle $\tilde E$ over $X_{dR}$.
\end{example} 
\begin{example}
    A transversal foliation on a scheme $Y$ over $X$ is the same as an isomorphism of $Y \simeq X \times _{X_{dR}}\tilde Y$ for a sheaf $\tilde Y$ over $X_{dR}$.
\end{example}

\begin{example}[Principal bundles]
    Let $G$ be a smooth affine algebraic group. Integrable connections on a principal $G$-bundle $P$ is equivalent to a principal $G$-bundle $\widetilde{P}$ on $X_{dR}$. To see this we just notice that the isomorphism $p_0^*P \simeq p_1^*P$ should be $G$-equivariant. It is classical that such connection is also equivalent to a $G$-equivariant $ad(P)$ valued $1$-form.
\end{example}

The Dolbeault stack $X_{Dol}$ is defined as $[X/\widehat{T_X}]$ where $T_X$ is the tangent bundle $\mathbf{Spec} (\operatorname{Sym \Omega_X})$, regardes as a group scheme over $X$, and $\widehat{TX}$ is the complextion of $T_X$ along the zero section. It may equally be regarded as the Cartier dual of $T^*X$, and we thus have
\[
D_{qc}(X_{Dol}) \simeq D_{qc}(T^*X),
\]
and therefore vector bundles on $X_{Dol}$ are the same as Higgs bundles.

\begin{example}\label{inthiggssheaf}
    Let us illustrates quasicoherent sheaves on $X_{Dol}$ are indeed Higgs bundles. This is quite similar to how to regard integrable connections as quasicoherent sheaves with stratifications, but even simpler. Given a Higgs bundle $(E,\theta)$ on $X$, we define a morphism
    \[
    \hat\theta \colon E \to E \otimes T_X^\wedge,\quad s \mapsto \sum_I \theta_I s \otimes dx^I
    \]
    where $x_1,\dots,x_n$ are etale coordinates of $X$. By extension of scalars we get the isomorphism $T_X^\wedge \otimes E \to E \otimes T_X^\wedge$.

    Conversely, if we have an isomorphism  $T_X^\wedge \otimes E \to E \otimes T_X^\wedge$ satisfying the cocycle condition which means , after resticrion of scalars,
    we have the commutative diagram
\[\begin{tikzcd}
	E & {E \otimes T_X^\wedge} \\
	{E \otimes T_X^\wedge} & {E \otimes T_X^\wedge \otimes T_X^\wedge}
	\arrow["{\hat\theta}", from=1-1, to=1-2]
	\arrow["{\hat \theta}"', from=1-1, to=2-1]
	\arrow["{\mathrm{id_E} \otimes m}", from=1-2, to=2-2]
	\arrow["{\hat \theta \otimes \mathrm{id}_{T_X^\wedge}}", from=2-1, to=2-2]
\end{tikzcd}\]
where $m$ is the multiplication on the group scheme $T_X^\wedge$. We take truncations 
\[\begin{tikzcd}
	E & {E \otimes T_X^{(2)}} \\
	{E \otimes T_X^{(1)}} & {E \otimes T_X^{(1)} \otimes T_X^{(1)}}
	\arrow["{\theta^{(2)}}", from=1-1, to=1-2]
	\arrow["{\theta^{(1)}}"', from=1-1, to=2-1]
	\arrow["{\mathrm{id_E} \otimes m}", from=1-2, to=2-2]
	\arrow["{ \theta^{(1)} \otimes \mathrm{id}_{T_X^{(1)}}}", from=2-1, to=2-2]
\end{tikzcd}\]
The map $\theta^{(1)}$ gives a homomorphism from $\theta \colon E \to E \otimes \Omega_X^\vee$. Take local coordinates $x_1,\dots,x_n$ of $X$, then $m(dx_i) = dx_i \otimes 1 + 1 \otimes dx_i$. Denote
\[
\theta^{(2)}(s) = s \otimes 1 + \theta_is \otimes dx^i + \Theta_{i,j}s \otimes dx^idx^j.
\]
Then the commutativity of the diagram tells us
\[
\theta_i\theta_j = \Theta_{i,j} = \theta_j\theta_i,
\]
and this is the integrability of Higgs field.
\end{example}
\begin{example}[Principal bundles]
    Let $G$ be a smooth affine algebraic group. Recall the definition of a Higgs field on a principal $G$-bundle $P$ is a Higgs field $\theta$ on $ad(P)$. Similar to the connection case, this is equivalent to a principal $G$-bundle $\widetilde{P}$ on $X_{Dol}$.
\end{example}

The Hodge stack is a deformation from the de Rham stack to the Dolbeault stack. It is defined by the groupoid object
\[\begin{tikzcd}
	R & {X \times \mathbb{A}^1}
	\arrow[shift left, from=1-1, to=1-2]
	\arrow[shift right, from=1-1, to=1-2]
\end{tikzcd}\]
where $R$ is the deformation to the normal cone of $X \hookrightarrow X \times X$, completed along $X \times \mathbb{A}^1$, and the composition is again induced by the tautological composition $p_{02}$.
The stack $X_{Hod}$ is a $\mathbb{G}_m$-equivariant stack over $\mathbb{A}^1$. Quasicoherent sheaves on $X_{Hod}$ are the same as integrable $t$-connections on $X \times \mathbb{A}^1$, and $\mathbb{G}_m$-equivariant sheaves on $X_{Hod}$ are the same as integrable connections on $X$, together with a Griffiths transverse filtration.

Now let us explain the relative case. Let $f \colon X \to S$ be a smooth surjection, then we can also form the stacks $(X/S)_{dR}, (X/S)_{Dol}, (X/S)_{Hod}$. For example, $(X/S)_{Hod}$ is defined as the pullback
\[\begin{tikzcd}
	{(X/S)_{Hod}} & {X_{Hod}} \\
	{S \times \mathbb{A}^1} & {S_{Hod}.}
	\arrow[from=1-1, to=1-2]
	\arrow[from=1-1, to=2-1]
	\arrow[from=1-2, to=2-2]
	\arrow[from=2-1, to=2-2]
\end{tikzcd}\]
Since $f \colon X \to S$ is smooth surjective, we may also define the relative stacks via pullback the presentation of groupoids. Let $T$ be an $S$-scheme, and $\lambda \in \Gamma(\mathcal{O}_T)$ defining $\lambda \colon T \to \mathbb{A}^1$, we may define $(X_T/T)_{\lambda,Hod}$ as
\[\begin{tikzcd}
	{(X_T/T)_{\lambda,Hod}} & {(X/S)_{Hod}} & {X_{Hod}} \\
	T & {S \times \mathbb{A}^1} & {S_{Hod}}
	\arrow[from=1-1, to=1-2]
	\arrow[from=1-1, to=2-1]
	\arrow[from=1-2, to=1-3]
	\arrow[from=1-2, to=2-2]
	\arrow[from=1-3, to=2-3]
	\arrow["{(\mathrm{id},\lambda)}", from=2-1, to=2-2]
	\arrow[from=2-2, to=2-3]
\end{tikzcd}\]

Quasicoherent sheaves on $(X_T/T)_{\lambda,Hod}$ are the same as integrable $\lambda$-connections on $X_T$, and $(X/S)_{1,Hod} = (X/S)_{dR}$, $(X/S)_{0,Hod} = (X/S)_{Dol}$.

If $f_0,f_1$ are two morphisms from $T$ to $S$, and $(f_0,f_1) \colon T \to S \times S$ factors through $D(1)$, then two maps $(f_0,\lambda)$ and $(f_\lambda,\lambda)$ from $T$ to $S_{Hod}$ are isomorphic, with the ismorphism given by $f_1$. And thus the pullback stacks $(f_0^*X/T)_{\lambda, Hod}$ and $(f_\lambda^*X/T)_{\lambda, Hod}$ are isomorphic. This should be the ultimate reason for the construction of the nonabelian $
\lambda$-connection in previous sections.

Let $G$ be a group scheme, $f \colon X \to S$ be morphism of schemes, and consider the induced map $f_{Hod} \colon X_{Hod} \to S_{Hod}$. Then $\mathcal{F} :=f_{Hod,*}(BG \times S_{Hod})$ is a $\mathbb{G}_m$-equivaraint sheaf over $S_{Hod}$, with the following properties:
\begin{itemize}
    \item The pullback of $\mathcal{F}$ to $S$ is the Hodge moduli stack $\mathcal{M}_{Hod}(X/S,G)$.
    \item The pullback of $\mathcal{F}$ to $S_{dR}$ is a sheaf on $S_{dR}$, whose further pullback to $S$ is $\mathcal{M}_{dR}(X/S,G)$.
    \item The pullback of $\mathcal{F}$ to $S_{Dol}$ is a sheaf on $S_{Dol}$, whose further pullback to $S$ is $\mathcal{M}_{Dol}(X/S,G)$
    \end{itemize}
 In this way, we can say we get a nonabelian Gauss-Manin connection on $\mathcal{M}_{dR}(X/S,G)$, together with a Griffiths transverse filtration, such that the graded object is $\mathcal{M}_{Dol}(X/S,G)$.

\begin{example}
    Let $f \colon X \to S$ be a smooth projective morphism. Then $R^kf_{Hod,*}\mathcal{O}_{X_{Hod}}$
    is the $k$-th algebraic de Rham cohomology of $X/S$, together with the Gauss--Manin connection and the Hodge filtration.
\end{example}

\begin{example}
    In the above example, we let $k = 1$. Then $R^1f_{Hod,*}\mathbb{G}_a \simeq \pi_0(f_{Hod,*}B\mathbb{G}_a)$
    is an isomorphism between the usual abelian Gauss--Manin connection and Hodge filtration to the non-abelian Gauss--Manin connection and Hodge filtration.
\end{example}

Let us take a closer look at the case $\lambda = 1,0$ to make things clearer. 
\begin{example}[Non-abelian Gauss--Manin connection]
    Let $X \to S$ be a smooth morphism, and $T$ be an $S$-scheme, then as the formation of de Rham stack commutes with base change, we have 
    \[
    X_{dR} \times _{S_{dR}} T \simeq (X_T/T)_{dR}.
    \]
    If $f_0,f_1$ are two morphisms from $T$ to $S$, such that $(f_0,f_1) \colon T \to S \times S$ factors through $(S \times S)^\wedge$, then $f_0,f_1$ defines the same morphism $F$ from $T \to S_{dR}$, whence we have
    \[
    (f_0^*X/T)_{dR} \simeq f_0^*(X/S)_{dR} \simeq F^*X_{dR} \simeq f_1^*(X/S)_{dR} \simeq (f_1^*X/T)_{dR}
    \]
    and thus we have an equivalence from the connections on $f_0^*X/T$ and the connections on $f_1^*X/T$. This is the non-abelian Gauss--Manin connection.
\end{example}
\begin{example}[Invariance of de Rham stacks under thickening]
    Let $S \to \tilde S$ be a thickening, and $\tilde X$ be a smooth lifting of $X$ to $\tilde S$. Then $(\tilde X/\tilde S)_{dR}$ is independent of the lifting $X$.
    In fact, the thickening $S \to \tilde S$ defines a morphism $\tilde S \to S_{dR}$, and we have
    \[
    (\tilde X/\tilde S)_{dR} \simeq X_{dR} \times_{S_{dR}} \tilde S.
    \]
    In particular, the category of integrable connections on $\tilde X/\tilde S$ is independent of the lifting $\tilde X$.
\end{example}
\begin{proposition}
    The invariance of integrable connections under thickening given in Proposition \ref{invariancefirstorder} is induced by the invariance of de Rham stack.
\end{proposition}
\begin{proof}
    This is the classical theory of crystals, see for example \cite{PMIHES_1994__80__5_0},Section 8.
\end{proof}

\begin{example}[Non-abelian Kodaira--Spencer]
    Let $X \to S$ be a smooth morphism, then we have Caetesian squares
\[\label{dolbeault}\tag{1}
\begin{tikzcd}
	{[X/T_{X/S}^\wedge]} & {[X/T_{X}^\wedge]} \\
	X & {[X/f^*T_S^\wedge]} \\
	S & {[S/T_S^\wedge]}
	\arrow[from=1-1, to=1-2]
	\arrow[from=1-1, to=2-1]
	\arrow[from=1-2, to=2-2]
	\arrow[from=2-1, to=2-2]
	\arrow[from=2-1, to=3-1]
	\arrow[from=2-2, to=3-2]
	\arrow[from=3-1, to=3-2]
\end{tikzcd}\]
So the non-abelian Kodaira Spencer map is induced from the action of $f^*T_S^\wedge$ on $[X/T_{X/S}^\wedge]$. From the Cartesian diagrams
\[\begin{tikzcd}
	{f^*T_S^\wedge} & X \\
	{[X/T_{X/S}^\wedge]} & {[X/T_X^\wedge]} \\
	X & {[X/f^*T_S^\wedge]}
	\arrow[from=1-1, to=1-2]
	\arrow[from=1-1, to=2-1]
	\arrow[from=1-2, to=2-2]
	\arrow[from=2-1, to=2-2]
	\arrow[from=2-1, to=3-1]
	\arrow[from=2-2, to=3-2]
	\arrow[from=3-1, to=3-2]
\end{tikzcd}\]
we know the $f^*T_S^\wedge$ action on $[X/T_{X/S}^\wedge]$ is induced by the group homomorphism $f^*T_S^\wedge \to [X/T_{X/S}^\wedge]$. If we omit the completion, then this is exactly the Kodaira--Spencer map.

Let us explain the action more explicitly. If $\tilde X$ is a stack over $BG$ and $X = *\times_{BG} \tilde X$, then the $G$-action on $X$ is defined via the dashed arrow in the cubical Cartesian diagram
\[\begin{tikzcd}
	{G \times X} && X \\
	& X && {\tilde{X}} \\
	G && {*} \\
	& {*} && BG
	\arrow[from=4-2, to=4-4]
	\arrow[from=2-4, to=4-4]
	\arrow[from=3-3, to=4-4]
	\arrow[from=3-1, to=3-3]
	\arrow[from=3-1, to=4-2]
	\arrow[from=1-3, to=2-4]
	\arrow[from=1-3, to=3-3]
	\arrow[from=1-1, to=1-3]
	\arrow[from=1-1, to=3-1]
	\arrow[from=2-2, to=2-4]
	\arrow[from=2-2, to=4-2]
	\arrow[dashed, from=1-1, to=2-2]
\end{tikzcd}\]
Let $H$ be a normal subgroup of $G$, then the pullback diagram
\[\begin{tikzcd}
	BH & BG \\
	{*} & BQ
	\arrow[from=1-1, to=1-2]
	\arrow[from=1-1, to=2-1]
	\arrow[from=2-1, to=2-2]
	\arrow[from=1-2, to=2-2]
\end{tikzcd}\]
where $Q:= G/H$ tells us there is a $Q$-action on $BH$. For any test scheme $T$, this Cartesian square identifies an $H$-torsor $Q$ over $T$ as a $G$-torsor $P$ over $T$ together with an isomorphism $\gamma \colon P/H \simeq T \times G/H$. Any $g \in G/H(T)$ acts on $\gamma$ via left multiplication. The difference of $H$-torsors given by $g \circ \gamma$ and $\gamma$ is the image of $g$ from $G/H(T) \to BH(T)$.
\end{example}
\begin{proposition}
    The auto-equivalence of the category of Higgs bundles on $X_0$ induced by $X_1$ as in \ref{invariancefirstorder} is induced by the $T_S^\wedge$-action on $(X/S)_{Dol}$.
\end{proposition}
\begin{proof}
    In general, let $G$ be a sheaf of commutative groups on $(Sch/X)_{et}$, then an element $v$ in $H^1(X,G)$ will induces an autoequivalence of $Qcoh(BG)$, since it will induce an automorphism of $BG$
    if we regard $v$ as an element inside $\pi_0(Map(*,BG))$. If we write $v$ as Cech cocycles $v_{ij}$, then the action on a quasicoherent sheaf $E$ on $BG$ is locally given by twisting via the coaction $V \to V \otimes_{\mathcal{O}_X}\mathcal{O}_G$ and sections $v_{ij}$.
\end{proof}

\section{Griffiths transversality}

Consider the open submoduli $M^{sm}_{dR}:=M^{sm}_{dR}(X/S,G) \to S$. Then the non-abelian Gauss--Manin connection will give a section of $\mathbf{T}_{M^{sm}_{dR}} \to \mathbf{T}_S$. This is described as follows (see also the appendix of \cite{lam2025algebraicityintegralitysolutionsdifferential}). Let $v \in \mathbf{T}_{S}$ be a tangent vector to $s \in S$. Let $X_v$ be the base change
\[\begin{tikzcd}
	{X_v} & X \\
	{k[\epsilon]} & S.
	\arrow[from=1-1, to=1-2]
	\arrow[from=1-1, to=2-1]
	\arrow[from=1-2, to=2-2]
	\arrow["v", from=2-1, to=2-2]
\end{tikzcd}\]
Let $(V,\nabla$) be a point on $M_{dR}$ lying above $s$, which means an integrable connection on $X_s$. Liftings of $v$ to $\mathbf{T}_{M_{dR}}(V,\nabla)$ are by definition integrable connections on $X_v$ whose reduction to $X_s$ are $(V,\nabla)$. We want to specify a distinguished element among them.

Indeed, there is a natural equivalence 
\[
\mathrm{Conn} (X[\epsilon]) \simeq \mathrm{Conn}(X_v)
\]
where $X[\epsilon]$ is the trivial deformation to $k[\epsilon]$, since both of them are flat liftings of $X_s$. And there is a distinguished element in $\mathrm{Conn} (X[\epsilon])$, namely the pull back of $(V,\nabla)$. This is the lifting of tangent vectors associated to non-abelian Gauss--Manin connection.
\begin{remark}
    Our construction gives a crystal structure on the moduli stack, or moduli functor. But since the good moduli space universally corepresents the moduli functor, the crystal structure will descend down to the good moduli space.
\end{remark}

\begin{proposition}
    The nonabelian Gauss--Manin connection is integrable. That means, $\nabla_{GM} \colon T_S \to T_{M_{dR}}$ preserves the Lie brackets, or equivalently $\operatorname{Im}(\nabla_{GM}) \subseteq T_{M_{dR}}$ is closed under the Lie bracket.
\end{proposition}
\begin{proof}
    This is because $M_{dR}$ underlies a crystal of schemes over $S$.
\end{proof}

As $\operatorname{Map}_{S_{Hod}}(X_{Hod}, BG) \to S_{Hod} $ is $\mathbb{G}_m$-equivariant, we may say  the Hodge filtration on $\mathcal{M}_{dR}(X/S,G)$ is Griffiths transverse with respect to the nonabelian Gauss--Manin connection. This is justified by the fact that $\mathbb{G}_m$-equivariant quasi-coherent sheaves on $S_{Hod}$ are the same as objects in $MIC(S)$ together with Griffiths transverse filtrations.

\begin{proposition}
    Suppose $X$ is defined over $S \times \mathbb{A}^1$, together with a $\mathbb{G}_m$-equivariant descent to $S_{Hod}$. This means there is a stack $\tilde{X}$ over $S_{Hod}$ such that there is a Cartesian diagram
\[\begin{tikzcd}
	X & {\tilde{X}} \\
	{S \times \mathbb{A}^1} & {S_{Hod}},
	\arrow[from=1-1, to=1-2]
	\arrow[from=1-1, to=2-1]
	\arrow[from=1-2, to=2-2]
	\arrow[from=2-1, to=2-2]
\end{tikzcd}\]
    and $\tilde X \to S_{Hod}$ is $\mathbb{G}_m$-equivariant.
    Then there is a $\mathbb{G}_m$-equivariant $t$-connection on $X$ over $S \times \mathbb{A}^1$.
\end{proposition}  
\begin{proof}
    Recall that $S_{Hod}$ is given by the quotient
\[\begin{tikzcd}
	R & {S \times \mathbb{A}^1}
	\arrow["t"', shift right, from=1-1, to=1-2]
	\arrow["s", shift left, from=1-1, to=1-2]
\end{tikzcd}\]
where $R$ is the deformation to the normal cone of $S \hookrightarrow S \times S$ completed along $S \times \{0\}$. If $X$ admits a $\mathbb{G}_m$-equivariant descent to $S_{Hod}$, then the morphism $s^*X \to t^*X$ is also $\mathbb{G}_m$-equivariant. To determine the lift of a tangent vector field, it suffices to work with the first order relation $R^1 \simeq D(1) \times \mathbb{A}^1$. Notice that we need to endow $D(1) \times \mathbb{A}^1$ a $\mathbb{G}_m$ action which is different from the obvious one, via the isomorphism to Rees deformation algebra.

Now let $\delta \colon Z[\epsilon] \to S \times \mathbb{A}^1$ be a tangent vector in $T_{S \times \mathbb{A}^1/\mathbb{A}^1}$, and let $\iota$ be the constant map associated to $\delta$. Then we have an induced map $(\iota, \delta) \colon Z[\epsilon] \to D(1) \times \mathbb{A}^1$. Recall the lifting of a tangent vector field is defined via
\[\begin{tikzcd}
	&& {p_t^*X} & {p_t^*X} \\
	& {p_0^*X} & {p_0^*X} \\
	{Z[\epsilon]} & {R^1} & {R^1}
	\arrow["{\lambda_*}", from=1-3, to=1-4]
	\arrow[from=2-2, to=1-3]
	\arrow["{\lambda_*}", from=2-2, to=2-3]
	\arrow[from=2-2, to=3-2]
	\arrow[from=2-3, to=1-4]
	\arrow[from=2-3, to=3-3]
	\arrow[from=3-1, to=2-2]
	\arrow[from=3-1, to=2-3]
	\arrow["{(\iota,\delta)}"', from=3-1, to=3-2]
	\arrow["{\lambda_*}"', from=3-2, to=3-3]
\end{tikzcd}\]
Claim that $\lambda_*(\iota,\delta) = (\lambda_*\iota, {\lambda}^{-1}\lambda_*\delta)$.
To see this, we assume $S = \operatorname{Spec} A$, $Z = \operatorname{Spec}B$. Let $r \colon P^1 \to A$ be the map induced by the diagonal embedding, and $\delta = f + \epsilon D$ from $A$ to $B [\epsilon]$, which induces a map $F \colon P^1 \to B[\epsilon]$, and $t \mapsto b \in B$. We have the isomorphism

\begin{align*}
    P^1[t] &\longrightarrow \Omega^1t^{-1} \oplus P^1[t],\\
    \sum a_k t^k &\longmapsto \sum (a_{k+1}-r(a_{k+1}) + r(a_k))t^k
\end{align*}
with the inverse given by
\[
\sum a_k't^k \longmapsto \sum (r(a_k')+a_{k-1}'-r(a_{k-1}'))t^k.
\]
Now $\lambda^*(\sum a_k' t^k) = \sum a_k'\lambda^kt^k$, and therefore 
\[
\lambda^* (\sum_k a_kt^k) = \sum_k (\lambda^kr(a_k') +\lambda^{k-1}(a_{k-1}' - r(a_{k-1}')))t^k.
\]
To summarize, the map $\lambda_*(\iota,\delta)$ is given by
\[
\sum_{k} a_k t^k \longmapsto \sum_{k} \lambda^k Fr(a_k') b^k + \sum_k \lambda^{k-1}Fr(a_{k-1}'-r(a_{k-1}'))b^k.
\]
But this is precisely the map induced by $(\lambda_*\iota, \lambda^{-1}\lambda_*\delta)$ (changing $b$ to $\lambda b$ and $D$ to $\lambda^{-1}D$).
\end{proof}

\section{The proof}
In the final section, we shall prove Theorem \ref{main result}. From discussions in the previous section, we know that 
$$
\theta_{KS}: g^*T_S\to T_{M^{sm}_{Dol}(X/S,G)/S} 
$$ 
is nothing but the nonabelian 0-connection over $M^{sm}_{Dol}:=M^{sm}_{Dol}(X/S,G)\to S$. We will use $\theta_{KS}$ for both the Higgs field and the associated vector field. It is described as follows. Let $v \in \mathbf{T}_{S}$ be a tangent vector to $s \in S$. Let $X_v$ be the base change as above. Let $(E,\theta)$ be a point in $M^{sm}_{Dol}$ lying over $s$, which means a Higgs bundle on $X_s$. Vectors in $\mathbf{T}_{M^{sm}_{Dol}/S}(E,\theta)$ are by definition Higgs bundles on $X[\epsilon]$ whose reduction to $X_s$ are $(E,\theta)$. 

The deformation $X_v$ gives an element in $H^1(X_s,T_{X_s})$, say given by Cech cocycle $\chi_{\alpha\beta} \in T_{X,\alpha\beta}:=T_{U_{\alpha\beta}}$. Then $\theta_{\chi_{\alpha\beta}}$ is an endomorphism of $E$, and $1 + \epsilon \theta_{\chi_{\alpha\beta}}$ is an automorphism of $E[\epsilon]$. Notice that $1 + \epsilon \chi_{\alpha\beta}$ is an automorphism of $U_{\alpha\beta}[\epsilon]$ which differs from $1$ by a square-zero ideal. From Proposition \ref{pullbackconnection} we get an automorphism
\[
    (E,\theta)[\epsilon]_{\alpha\beta} \simeq (E,\theta)[\epsilon]_{\alpha\beta}, \quad s \mapsto s + \sum_{i} \theta_i(s)\epsilon \chi_{\alpha\beta}(x_i)
\]
However, $
\chi_{\alpha\beta}(x_i)\frac{\partial}{\partial x_i} = \chi_{\alpha\beta}$, and therefore this automorphism is just $1 + \epsilon \theta_{\chi_{\alpha\beta}}$.

By abuse of notation in the following we also use $X$ to denote $X_s$. We have constructed a map from deformations of $X$ to deformations of $(E,\theta)$. Recall the deformations of Higgs bundle are controlled by $\mathbb{H}^1(\Omega^*_{Hig}(ad(E),\theta^{ad}))$ where
\[\begin{tikzcd}
	{ad(E)} & {ad(E) \otimes \Omega^1} & {ad(E)\otimes \Omega^2}.
	\arrow["{\theta^{ad}}", from=1-1, to=1-2]
	\arrow["{\theta^{ad}}", from=1-2, to=1-3]
\end{tikzcd}\]
Take an affine open cover $U_{\alpha}$ of $X$, and consider the Cech resolution
\[\begin{tikzcd}
	{ad(E)} & {ad(E) \otimes \Omega^1} & {ad(E) \otimes \Omega^2} \\
	{ad(E)_{\alpha}} & {{ad(E) \otimes \Omega^1}_{\alpha}} & {{ad(E) \otimes \Omega^2}_{\alpha}} \\
	{ad(E)_{\alpha\beta}} & {{ad(E) \otimes \Omega^1}_{\alpha\beta}} \\
	{ad(E)_{\alpha\beta\gamma}}
	\arrow[from=1-1, to=1-2]
	\arrow[from=1-1, to=2-1]
	\arrow[from=1-2, to=1-3]
	\arrow[from=1-2, to=2-2]
	\arrow[from=1-3, to=2-3]
	\arrow[from=2-1, to=2-2]
	\arrow[from=2-1, to=3-1]
	\arrow[from=2-2, to=2-3]
	\arrow[from=2-2, to=3-2]
	\arrow[from=3-1, to=3-2]
	\arrow[from=3-1, to=4-1]
\end{tikzcd}\]
An element in $\mathbb{H}^1$ is represented by $s_{\alpha\beta} \in ad(E)_{\alpha\beta}$ and $t_\alpha \in {ad(E) \otimes \Omega^1}_{\alpha}$ such that
\begin{enumerate}
    \item $s_{\alpha\beta} + s_{\beta\gamma} + s_{\gamma\alpha} = 0$,
    \item $t_\alpha - t_\beta = m_\theta (s_{\alpha\beta})$,
    \item $\theta^{ad}(t_\alpha) = 0$.
\end{enumerate}
Given such $s_{\alpha\beta},t_\alpha$ the deformation of Higgs bundle is constructed as follows:
The morphism $1 + \epsilon s_{\alpha\beta}$ is an automorphism of $E[\epsilon]_{\alpha\beta}$, and $\theta_\alpha + \epsilon t_\alpha$ defines a Higgs field on $E$, compatible with automorphisms $1 + \epsilon s_{\alpha\beta}$. These glue together to a Higgs bundle on $X[\epsilon]$.

There is a natural map $\theta \colon T_{X_s} \to ad(E)$, which induces a map 
$$
H^1(X_s,T_{X_s}) \to \mathbb{H}^1(\Omega^*_{Hig}(ad(E),\theta^{ad})).
$$
\begin{proposition}
    $\tau_s$,  the value of the morphism $\tau$ (see \S1) at $s$, is the map 
    $$H^1(X_s,T_{X_s}) \to \mathbb{H}^1(\Omega^*_{Hig}(ad(E),\theta^{ad})).$$
\end{proposition}
\begin{proof}
\[\begin{tikzcd}
	{T_X} & {ad(E)} & {ad(E) \otimes \Omega^1} & {ad(E) \otimes \Omega^2} \\
	{T_{X,\alpha}} & {ad(E)_{\alpha}} & {{ad(E) \otimes \Omega^1}_{\alpha}} & {{ad(E) \otimes \Omega^2}_{\alpha}} \\
	{T_{X,\alpha\beta}} & {ad(E)_{\alpha\beta}} & {{ad(E) \otimes \Omega^1}_{\alpha\beta}} \\
	{T_{X,\alpha\beta\gamma}} & {ad(E)_{\alpha\beta\gamma}}
	\arrow[from=1-1, to=1-2]
	\arrow[from=1-1, to=2-1]
	\arrow[from=1-2, to=1-3]
	\arrow[from=1-2, to=2-2]
	\arrow[from=1-3, to=1-4]
	\arrow[from=1-3, to=2-3]
	\arrow[from=1-4, to=2-4]
	\arrow[from=2-1, to=2-2]
	\arrow[from=2-1, to=3-1]
	\arrow[from=2-2, to=2-3]
	\arrow[from=2-2, to=3-2]
	\arrow[from=2-3, to=2-4]
	\arrow[from=2-3, to=3-3]
	\arrow[from=3-1, to=3-2]
	\arrow[from=3-1, to=4-1]
	\arrow[from=3-2, to=3-3]
	\arrow[from=3-2, to=4-2]
	\arrow[from=4-1, to=4-2]
\end{tikzcd}\]
Let $\chi_{\alpha\beta}$ be a cocycle representing an element in $H^1(X,T_X)$. Then the corresponding $s_{\alpha\beta} = \theta_{\chi_{\alpha\beta}}$ and $t_\alpha = 0$. Thus the deformed Higgs bundle is $(E,\theta)[\epsilon]_{\alpha}$ glued via automorphisms $1 + \epsilon \theta_{\chi_{\alpha\beta}}$, which is exactly given by the non-abelian Kodaira-Spencer map.
\end{proof}

\begin{remark}
    This proof actually shows the formula holds on the stable locus.
\end{remark}

\begin{proposition} The non-abelian Kodaira Spencer map $\theta_{KS}$ is integrable.
\end{proposition}
\begin{proof}
    We know $\mathcal{M}_{Dol}(X/S,G)$ is isomorphic to $\widetilde{\mathcal{M}} \times_{S_{Dol}} S$ where 
    \[
    \widetilde{\mathcal{M}} = Map_{S_{Dol}}(X_{Dol},BG \times S_{Dol}).
    \]
    Since the moduli space is good, there is a sheaf $\widetilde M$ over $S_{Dol}$ such that 
    $$M_{Dol}(X/S,G) \simeq \widetilde{M} \times_{S_{Dol}}S.$$ Equivalently, there is an isomorphism
    \[
    M_{Dol}(X/S,G) \times_S TS^\wedge \simeq TS^\wedge \times_S M_{Dol}(X/S,G)
    \]
    satisfying the cocycle condition. Following the discussion made in \ref{inthiggssheaf}, this implies the integrability claimed in Theorem \ref{main result}.
\end{proof}

\begin{proposition}
    $\theta_{KS}$ is a graded Higgs field; that is, $(M^{sm}_{Dol},\theta_{KS}) \simeq (M^{sm}_{Dol},t\theta_{KS})$ for any $t \in k^*$. Indeed, we have the following commuatative diagram
\[\begin{tikzcd}
	{\mathcal{O}_{M^{sm}_{Dol}}} & {\mathcal{O}_{M^{sm}_{Dol}} \otimes \Omega_{S/k}^1} \\
	{\mathcal{O}_{M^{sm}_{Dol}}} & {\mathcal{O}_{M^{sm}_{Dol}} \otimes \Omega_{S/k}^1.}
	\arrow["{\theta_{KS}}", from=1-1, to=1-2]
	\arrow["t^*"', from=1-1, to=2-1]
	\arrow["{t^* \otimes \mathrm{id}}", from=1-2, to=2-2]
	\arrow["{t\theta_{KS}}", from=2-1, to=2-2]
\end{tikzcd}\]
Or in terms of associated lifting of vector fields, we have the following commutative diagram
\[\begin{tikzcd}
	{g^*T_S} & {t^*g^*T_S} \\
	{T_{M_{Dol}/S}} & {t^*T_{M_{Dol}/S}.}
	\arrow[from=1-1, to=1-2]
	\arrow["{t\theta_{KS}}"', from=1-1, to=2-1]
	\arrow["{t^*\theta_{KS}}", from=1-2, to=2-2]
	\arrow["{t_*}", from=2-1, to=2-2]
\end{tikzcd}\]
\end{proposition}

\begin{proof} This already follows from proposition \ref{gradedness} and the Griffiths transversality. We also give a direct proof using the explicit formula here.
 Given $t \in k^*$, $(E,\theta)$ a Higgs bundle on $X_s$, we know $t(E,\theta) = (E,t\theta)$ and the induced map on tangent space
\[
\text{deformations of}\ (E,\theta) \to 
\text{deformations of}\ (E,t\theta)
\]
is $(\tilde E,\tilde \theta) \mapsto (\tilde E, t \tilde\theta)$. We need to check the following commutative diagram
\[\begin{tikzcd}
	{H^1(X_s,T_{X_s})} & {\mathbb{H}^1(ad(E)_{\theta-{Higgs}})} \\
	{H^1(X_s,T_{X_s})} & {\mathbb{H}^1(ad(E)_{t\theta-Higgs})}
	\arrow["{t\theta_{KS}}", from=1-1, to=1-2]
	\arrow[shift right, no head, from=1-1, to=2-1]
	\arrow[shift left, no head, from=1-1, to=2-1]
	\arrow["t", from=1-2, to=2-2]
	\arrow["{\theta_{KS}}", from=2-1, to=2-2]
\end{tikzcd}\]
But this follows from the explicit description of the non-abelian Higgs field, namely given $\chi_{\alpha\beta}\in H^1(X_s,T_{X_s})$, the deformation of $(E,\theta)$ is given by twisting of $1 + \epsilon \theta_{t\chi_{ij}}$ and the deformation of $(E,t\theta)$ is given by twisting of $1 + \epsilon t\theta_{\chi_{ij}}$.
\end{proof}

\begin{corollary}
    If $(E,\theta)$ is a graded Higgs bundle on $X_s$, then $\theta_{KS}|_{(E,\theta)}$ lies on the weight $-1$ subspace of the $\mathbb{G}_m$-representation $T_{(E,\theta)}M_{Dol}(X_s)$.
\end{corollary}

The $\mathbb{G}_m$-decomposition of the tangent space at $\mathbb{G}_m$-fixed points, namely graded Higgs bundles, can be described explicitly. Indeed, for any Higgs bundle $(E,\theta)$, the tangent differential $t_*$ on $\mathbb{H}^1$ is induced on the level of complex
\[\begin{tikzcd}
	{\operatorname{End}(E)} & {\operatorname{End}(E) \otimes\Omega^1} & {\operatorname{End}(E) \otimes \Omega^2} \\
	{\operatorname{End}(E)} & {\operatorname{End}(E)\otimes \Omega^1} & {\operatorname{End}(E) \otimes \Omega^2}
	\arrow["{ad_\theta}", from=1-1, to=1-2]
	\arrow["{\mathrm{id}}"', from=1-1, to=2-1]
	\arrow["{ad_\theta}", from=1-2, to=1-3]
	\arrow["{t }"', from=1-2, to=2-2]
	\arrow["{t^2 }"', from=1-3, to=2-3]
	\arrow["{ad_{t\theta}}", from=2-1, to=2-2]
	\arrow["{ad_{t\theta}}", from=2-2, to=2-3]
\end{tikzcd}\]
If $(E,\theta)$ is graded, that is, $E \simeq \bigoplus_{i = 0}^n E_i$, and $\theta \colon E_i \to E_{i+1} \otimes \Omega^1$. We let $\mathbb{G}_m$ act on $E$ as $\rho(t)|E_i = t^i \cdot \mathrm{id}_{E_i}$. Then there is a commutative diagram
\[\begin{tikzcd}
	E & {E \otimes \Omega^1} \\
	E & {E\otimes \Omega^1}
	\arrow["\theta", from=1-1, to=1-2]
	\arrow["{\rho(t)}"', from=1-1, to=2-1]
	\arrow["{\rho(t)}", from=1-2, to=2-2]
	\arrow["{t\theta}", from=2-1, to=2-2]
\end{tikzcd}\]
which shows $(E,\theta)$ and $(E,t\theta)$ are isomorphic Higgs bundles. Similarly, we can furthur get the diagram
\[\begin{tikzcd}
	{\operatorname{End}(E)} & {\operatorname{End}(E) \otimes\Omega^1} & {\operatorname{End}(E) \otimes \Omega^2} \\
	{\operatorname{End}(E)} & {\operatorname{End}(E)\otimes \Omega^1} & {\operatorname{End}(E) \otimes \Omega^2} \\
	{\operatorname{End}(E)} & {\operatorname{End}(E)\otimes \Omega^1} & {\operatorname{End}(E)\otimes \Omega^2}
	\arrow["{ad_\theta}", from=1-1, to=1-2]
	\arrow["{\mathrm{id}}"', from=1-1, to=2-1]
	\arrow["{ad_\theta}", from=1-2, to=1-3]
	\arrow["t"', from=1-2, to=2-2]
	\arrow["{t^2}"', from=1-3, to=2-3]
	\arrow["{ad_{t\theta}}", from=2-1, to=2-2]
	\arrow["{\rho(t^{-1})}"', from=2-1, to=3-1]
	\arrow["{ad_{t\theta}}", from=2-2, to=2-3]
	\arrow["{\rho(t^{-1})}"', from=2-2, to=3-2]
	\arrow["{\rho(t^{-1})}"', from=2-3, to=3-3]
	\arrow["{ad_\theta}", from=3-1, to=3-2]
	\arrow["{ad_\theta}", from=3-2, to=3-3]
\end{tikzcd}\]

The vertical composed morphism is just the $t_*$ action on $T_{(E,\theta)} M_{Dol}$. Notice that we need to use $(t^{-1})^*$ to identify $(E,t\theta)$ with $(E,\theta)$. Therefore we can decompose the deformation complex as 
\[
\operatorname{End}(E)_{Higgs} = \bigoplus_k \mathbf{U}_k^\bullet
\]
where 
\[
\mathbf{U}_k^\bullet := U_{k} \to U_{k+1} \otimes \Omega^1 \to U_{k+2} \otimes \Omega^2,
\]
and $U_k$ is the weight $k$ subspace of $\operatorname{End}(E) \simeq E^\vee \otimes E$. Due to the action of $t^{-1}$, the weight $-k$ subspace of $\mathbb{H}^1(\operatorname{End}(E)_{Higgs})$ is $\mathbb{H}^1(\mathbf{U}_k^\bullet)$.

\begin{example}
    Let $X$ be a smooth projective curve of genus $g \ge 2$. The uniformizing Higgs bundle $E$ is a graded Higgs bundle with $E \simeq K_X^{1/2} \oplus K_X^{-1/2}$, and $\theta \colon K_X^{1/2} \to K_X^{-1/2} \otimes K_X \simeq K_X^{1/2}$ is the identity. The Higgs bundle $\operatorname{Sym}^nE$ is a stable graded Higgs bundle of rank $n$ with underlying bundle
 $\bigoplus_{j=0}^n K_X^{\frac{n-2j}{2}}$, and $\theta \colon K_X^{\frac{n-2j}{2}} \to K_X^{\frac{n-2j-2}{2}} \otimes K_X$ is the identity.
 The weight $-1$ subspace is $\mathbb{H}^1(C^\bullet)$
 where 
 \[
 C^\bullet = \bigoplus_{j=0}^{n-1}K_X^{-1} \to \bigoplus_{j=0}^{n-2}K_X^{-1},\quad (x_j) \mapsto (x_j - x_{j+1}).
 \]
 The natural map $K_X^{-1} \to C_\bullet$ is a quasi-isomorphism, and thus $H^1(X,T_X) \simeq \mathbb{H}^1(C^\bullet)$.
\end{example}

\begin{proposition}
    Let $X \to S$ be a family of smooth projective curves of genus $g \ge 2$. If the nonabelian Kodaira Spencer $\theta_{KS} = 0$ on the moduli space $M_{Dol}(X/S,r)$ of Higgs bundles of rank $r \ge 2$, then the absolute nonabelian Kodaira-Spencer map $\rho_{KS} = 0$.
\end{proposition}
\begin{proof}
    For $s \in S$, consider the $F_s: = \operatorname{Sym}^{r-1}(E_s)$ where $E_s$ is the uniformizing Higgs bundle of $X_s$. The above calculation shows $H^1(X_s,T_{X_s}) \to \mathbb{H}^1(End(F)_{Higgs})$ is injective.
\end{proof}

\section{Appendix: Associated graded of abelian Gauss--Manin connection}
Let $f \colon X \to S$ be a smooth projective morphism, and we also use $f$ to denote the induced morphism on deRham, Dolbeault and Hodge stacks. Then $R^kf_* \mathcal{O}_{X_{Hod}}$ is the Rees module associated to the Gauss--Manin connection on $R^kf_*(\Omega_{X/S}^\bullet)$ together with the Hodge filtration. To compute the associated graded object, it suffices to compute the functor $R^kf_*$ for Dolbeault stacks. Recall the derived direct image of Higgs bundles defined in \cite[3.3]{PMIHES_2007__106__1_0}
\[
Rf_{Dol,*}: = R\bar f_{*} \circ Ri^! \colon D_{qcoh}(T^*X) \to D_{qcoh}(T^*S)
\]
where $i, \bar f$ are indicated in the following diagram
\[\begin{tikzcd}
	{T^*S \times_S X} & {T^*X} \\
	{T^*S}
	\arrow["i", hook, from=1-1, to=1-2]
	\arrow["{\bar f}"', from=1-1, to=2-1]
\end{tikzcd}\]
This can also be interpreted directly as $Rf_{Dol,*}$ where $f_{Dol,*} \colon X_{Dol} \to S_{Dol}$ is the induced morphism on Dolbeault stacks. To see this, we only need to consider the diagram
\[\begin{tikzcd}
	{X_{Dol}} & {[X/f^*T_S^\wedge]} \\
	& {S_{Dol}}
	\arrow["i^\vee", from=1-1, to=1-2]
	\arrow["\tilde{f}", from=1-2, to=2-2]
\end{tikzcd}\]
and notice that $Ri^! = Ri^\vee_*$, $R\bar f_* = R\tilde{f}_*$. Now it follows from \cite[Remark 3.3]{PMIHES_2007__106__1_0} that $R^kf_{*}\mathcal{O}_{X_{Dol}}$ is the $k$-the Hodge cohomology coupled with Kodaira--Spencer map.
\small

\bibliographystyle{unsrt}
\bibliography{reference}

@book {Voisin,
AUTHOR = {Voisin, Claire},
     TITLE = {Hodge theory and complex algebraic geometry. {I}},
    SERIES = {Cambridge Studies in Advanced Mathematics},
    VOLUME = {76},
   EDITION = {English},
      NOTE = {Translated from the French by Leila Schneps},
 PUBLISHER = {Cambridge University Press, Cambridge},
      YEAR = {2007},
     PAGES = {x+322},
      ISBN = {978-0-521-71801-1},
   MRCLASS = {32J25 (14C30 14D07 32G20)},
  MRNUMBER = {2451566},
}

@book{Simpson,
    author = {{Simpson}, Carlos},
    title = "{The Hodge filtration on nonabelian cohomology}",
    publisher = {Algebraic Geometry Santa Cruz 1995, Part 2. Proceedings of Symposia in Pure Mathematics},
    year = 1997,
}

@article{ChenTing,
author = {Chen, Ting},
year = {2012},
month = {05},
pages = {1-15},
title = {The associated map of the nonabelian Gauss-Manin connection},
volume = {10},
journal = {Central European Journal of Mathematics},
doi = {10.2478/s11533-011-0110-3}
}

@article {Katz,
    AUTHOR = {Katz, Nicholas M.},
     TITLE = {Algebraic solutions of differential equations ({$p$}-curvature
              and the {H}odge filtration)},
   JOURNAL = {Invent. Math.},
  FJOURNAL = {Inventiones Mathematicae},
    VOLUME = {18},
      YEAR = {1972},
     PAGES = {1--118},
      ISSN = {0020-9910,1432-1297},
   MRCLASS = {14D05 (14C30 14F30 34A20)},
  MRNUMBER = {337959},
MRREVIEWER = {T.\ Oda},
       DOI = {10.1007/BF01389714},
       URL = {https://doi.org/10.1007/BF01389714},
}

@misc{lam2025algebraicityintegralitysolutionsdifferential,
      title={Algebraicity and integrality of solutions to differential equations}, 
      author={Yeuk Hay Joshua Lam and Daniel Litt},
      year={2025},
      eprint={2501.13175},
      archivePrefix={arXiv},
      primaryClass={math.AG},
      note={\url{https://arxiv.org/abs/2501.13175}} 
}

@article{PMIHES_1994__80__5_0,
     author = {Simpson, Carlos T.},
     title = {Moduli of representations of the fundamental group of a smooth projective variety {II}},
     journal = {Publications Math\'ematiques de l'IH\'ES},
     pages = {5--79},
     publisher = {Institut des Hautes Etudes Scientifiques},
     volume = {80},
     year = {1994},
     mrnumber = {1320603},
     zbl = {0891.14006},
     language = {en},
     url = {https://www.numdam.org/item/PMIHES_1994__80__5_0/}
}

@article{PMIHES_2007__106__1_0,
     author = {Ogus, A. and Vologodsky, V.},
     title = {Nonabelian {Hodge} theory in characteristic $p$},
     journal = {Publications Math\'ematiques de l'IH\'ES},
     pages = {1--138},
     publisher = {Springer},
     volume = {106},
     year = {2007},
     doi = {10.1007/s10240-007-0010-z},
     language = {en},
     url = {https://www.numdam.org/articles/10.1007/s10240-007-0010-z/}
}

\end{document}